\documentclass[11pt,twoside,reqno]{amsart}
\usepackage{mathptmx,amsmath,amssymb,amsfonts,amsthm,
mathptmx,enumerate,color}
\usepackage{graphicx}

\usepackage{epstopdf}

\newtheorem{theorem}{Theorem}[section]
\newtheorem{corollary}{Corollary}[section]
\newtheorem{lemma}{Lemma}[section]

\theoremstyle{definition}

\newtheorem{remark}{Remark}[section]

\numberwithin{equation}{section}

\def\mychoose{\atopwithdelims[]}

\newcommand{\qbinomial}[3]{\mbox{$
\biggl[ 
\begin{array}{c}
#1\\
 #2
\end{array}\biggr]_{
\!{#3}}$}}

\newcommand{\be}{\begin{equation}}
\newcommand{\ee}{\end{equation}}
\newcommand{\bea}{\begin{eqnarray}}
\newcommand{\eea}{\end{eqnarray}}
\newcommand{\bd}{\begin{displaymath}}
\newcommand{\ed}{\end{displaymath}}

\begin{document}

\setcounter{page}{1}

\vspace*{1.0cm}

\title[$q$-difference equations  for   homogeneous  $q$-difference operators  and  their applications]{
$q$-difference equations  for  homogeneous  $q$-difference operators  and  their applications}

\author[  Sama Arjika]
{ Sama Arjika }

\maketitle

\vspace*{-0.6cm}

\begin{center}
{\footnotesize {\it
 Department of Mathematics and Informatics, 
University of Agadez,\break 
Post Box 199, Agadez, Niger}}
\end{center}

\vskip 4mm 

{\footnotesize \noindent {\bf Abstract.}
 
In this short paper, we show how to deduce several types of generating functions from  Srivastava  {\it et al} [Appl. Set-Valued Anal. Optim. {\bf 1} (2019), pp. 187-201.]  by the method of  $q$-difference equations. Moreover,  we build relations between transformation formulas and homogeneous $q$-difference equations. 

\vskip 1mm

\vskip 1mm 

}

\renewcommand{\thefootnote}{}
\footnotetext{ 
\par
E-mail addresses: rjksama2008@gmail.com (Sama Arjika). 
%
%
}

\section{Introduction and basics properties}

In this paper, we adopt the common conventions and notations on $q$-series.  For the convenience of the reader, we provide  a  summary of the mathematical notations, basics properties  and definitions to be used  in the sequel.   We refer  to the general references (see \cite{Koekock})  for the definitions and notations. Throughout this paper, we assume that  $|q|< 1$. 

For complex  numbers $a$, the $q$-shifted factorials are defined by:
\bea
  (a;q)_{n} =\left\{\begin{array}{ll}1& \quad \mbox{ if } n=0\\
 (1-a) (1-aq)\cdots (1-aq^{n-1}), & \quad \mbox{ if } n =1, 2, 3, \ldots\end{array}\right.
\eea
and  for tends  to infinity, we have  
$$
 (a;q)_{\infty}:=\prod_{k=0}^{\infty}(1-aq^{k}). 
$$
The following easily verified identities will be frequently used in this paper:
\be
 (a;q)_{n}=\frac{ (a;q)_{\infty}}{ (aq^n;q)_{\infty}}\quad 
 (a;q)_{n+k}=  (a;q)_n(aq^n;q)_k 
\ee
 and 
 $ (a_1,a_2, \ldots, a_r;q)_m=(a_1;q)_m (a_2;q)_m\cdots(a_r;q)_m,\; m\in\{0, 1, 2\cdots\}$.  \\
The $q$-binomial coefficients   are given by  
\bea 
\label{samz}
 {\,n\,\atopwithdelims []\,k\,}_{q} : =\left\{\begin{array}{ll}\frac{(q;q)_n }
{(q;q)_k\,(q;q)_{n-k}}& \mbox{ if } 0\leq k\leq n\\
0&\mbox{ otherwise}.\end{array}\right. \nonumber
\eea
  
The basic (or $q$-) hypergeometric function 
of the variable $z$ and with $\mathfrak{r}$ numerator 
and $\mathfrak{s}$ denominator parameters 
(see, for details, the monographs by 
Slater \cite[Chapter 3]{SLATER} 
and by Srivastava and Karlsson 
\cite[p. 347, Eq. (272)]{SrivastaKarlsson}; 
see also \cite{Koekock}) 
is defined as follows:
 $$
{}_{\mathfrak r}\Phi_{\mathfrak s}\left[\begin{array}{rr}a_1, a_2,\ldots, a_{\mathfrak r };\\
 \\
b_1,b_2,\ldots,b_{\mathfrak s};
 \end{array}\,
q;z\right]
 =\sum_{n=0}^\infty\Big[(-1)^n q^{({}^n_2)}\Big]^{1+{\mathfrak s}-{\mathfrak r }}\,\frac{(a_1, a_2,\ldots, a_{\mathfrak r};q)_n}{(b_1,b_2,\ldots,b_{\mathfrak s};q)_n}\frac{ z^n}{(q;q)_n}
$$
where $q\neq 0$ when ${\mathfrak r }>{\mathfrak s}+1$. Note that:
$$
{}_{\mathfrak r+1}\Phi_{\mathfrak r}\left[\begin{array}{rr}a_1, a_2,\ldots, a_{\mathfrak r+1}\\
 \\
b_1,b_2,\ldots,b_{\mathfrak r };
 \end{array}\,q;z\right]
 =\sum_{n=0}^\infty \frac{(a_1, a_2,\ldots, a_{\mathfrak r+1};q)_n}{(b_1,b_2,\ldots,b_{\mathfrak r};q)_n}\frac{ z^n}{(q;q)_n}.
$$
 
Here, in our present investigation, we are mainly concerned  with the Cauchy polynomials $p_n(x,y)$ as given below (see \cite{Chen2003,GasparRahman}): 
\bea
\label{def}
 p_n(x,y):=(x-y)( x- qy)\cdots ( x-q^{n-1}y) =(y/x;q)_nx^n
\eea
with the generating function \cite{Chen2003}
\be
\label{gener}
\sum_{n=0}^{\infty} p_n(x,y)
\frac{t^n }{(q;q)_n} = 
\frac{(yt;q)_\infty}{(xt;q)_\infty},
\ee
where \cite{Chen2003}
$$
  p_n(x,y)=(-1)^n q^{({}^n_2)}p_n(y,q^{1-n}x),
$$
and
$$
  p_{n-k}(x,q^{1-n}y)=(-1)^{n-k} q^{({}^k_2)-({}^n_2)}p_{n-k}(y,q^{k}x)
$$
 which naturally arise in the $q$-umbral calculus  \cite{Andrews}, Goldman and 
 Rota \cite{Goldman},  Ihrig and Ismail \cite{Ismail1981}, Johnson \cite{Johnson} 
 and Roman \cite{Roman1982}. The generating 
 function (\ref{gener}) is also the homogeneous version
of the Cauchy identity or the $q$-binomial theorem  \cite{GasparRahman}
\be
\label{putt}
\sum_{k=0}^{\infty} 
\frac{(a;q)_k }{(q;q)_k}z^{k}={}_{1}\Phi_0\left[\begin{array}{c}a
 \\
-
 \end{array};
q,z\right]= 
\frac{(az;q)_\infty}{(z;q)_\infty}\quad |z|<1. 
\ee
Putting    $a=0$, the relation (\ref{putt}) becomes  Euler's identity  \cite{GasparRahman}
\be
\label{q-expo-alpha}
  \sum_{k=0}^{\infty} \frac{ z^{k}}{(q;q)_k}=\frac{1}{(z;q)_\infty} \quad |z|<1
\ee
and its inverse relation  \cite{GasparRahman}
\be
\label{q-Expo-alpha}
 \sum_{k=0}^{\infty} 
\frac{(-1)^kq^{ ({}^k_2)
}\,z^{k}}{(q;q)_k}=(z;q)_\infty.
\ee

The following two $q$-difference operators are defined by \cite{Liu97,SrivastaAbdlhusein,Saadsukhi}
\be 
\label{deffd}
 {D}_q\big\{f(x)\big\}=\frac{f(x)-f( q x)}{ x}, \quad  
 {\theta}_{x}= {\theta_{xy}}_{|y=0},\quad {\theta}_{xy}\big\{f(x,y)\}:=\frac{f(q^{-1}x,y)-f( x,qy)}{q^{-1}x-y}.
\ee
The   Leibniz rule for the $D_q$ is the following identity\cite{Roman1982}
\be
\label{Lieb}
D_q^n\left\{ f(x)g(x)\right\}=\sum_{k=0}^n \qbinomial{n}{k}{q} q^{k(k-n)} D_q^k\left\{f(x)\right\} D_q^{n-k}\left\{g(q^kx)\right\}
\ee 
where $D_q^0$ is  understood as the identity.
For $f(x)=x^k$ and $g(x)=1/(x t;q)_\infty$, we have 
\be 
\label{tLieb}
D_q^n\left\{ \frac{x^k}{(x t;q)_\infty} \right\}= \frac{(q;q)_k }{(x t;q)_\infty}\sum_{j=0}^n \qbinomial{n}{j}{q}   \frac{ (x t;q)_j }{(q;q)_{k-j}}t^{n-j}x^{k-j}.
\ee 
  Saad and Sukhi  \cite{Saadsukhi,Saad} and Chen and Liu \cite{Liu97,Liu98}  employed the technique of parameter augmentation by constructing the following  $q$-exponential operators  
\be 
\label{mats}
R(bD_q)=\sum_{k=0}^\infty \frac{(-1)^kq^{({}^k_2)}  }{(q;q)_k}\, \left(b\,{D}_q\right)^k,\quad
\mathbb{E}(b\theta_{a})=\sum_{k=0}^\infty \frac{q^{({}^k_2) } }{(q;q)_k}\, \left(b\,\theta_{a}\right)^k.
\ee

\begin{theorem}(\cite[Theorem 2]{Liu10})
Let $f(a,b)$ be a two-variable analytic function in a neighbourhood of $(a,b)=(0,0)\in\mathbb{C}^2$.    If $f(a,b)$  satisfies the $q$-difference equation
\be 
\label{jjZA}
 a f(aq, b)-bf(a,  bq) 
 = (a-b) f(aq,  bq) 
\ee
then we have:
\bea
f(a,b)=\mathbb{E}(b\theta_{a})\Big\{f(a , 0) \Big\}.
\eea
\end{theorem}
 
 Liu \cite{Liu10,Liu11} initiated the method of $q$-difference equations  and 
deduced several results involving Bailey’s ${}_6\psi_6$, 
$q$-Mehler formulas for Rogers-Szeg\"o polynomials 
and $q$-integral of Sears’ transformation.

Recently,   Srivastava, Arjika and Kelil  \cite{ArjikaSri},  introduced  two homogeneous $q$-difference operators  $\widetilde{E}(a,b;D_{q})$ and $\widetilde{L}(a,b;\theta_{xy})$  
\be 
\label{aqamat}\widetilde{E}(a,b; {D}_{q})=\sum_{k=0}^\infty \frac{(-1)^kq^{({}^k_2)}\,(a;q)_k }{(q;q)_k}\, \left(b\,{D}_{q}\right)^k,\quad 
\widetilde{L}(a,b; \theta_{xy})=\sum_{k=0}^\infty \frac{q^{({}^k_2)}\,(a;q)_k }{(q;q)_k}\, \left(b\,\theta_{xy}\right)^k.
\ee
which turn out to be suitable for dealing with  a generalized Cauchy polynomials $p_n(x,y,a)$ \cite{ArjikaSri}
\be
\label{defR}
p_n(x,y,a) =\widetilde{E}(a,y; {D}_{q})\{ x^n\}.
\ee

The method of $q$-exponential operator is a rich and powerful tool for $q$-series, especially it makes many famous results easily fall into this framework. In this paper, we use this method to   derive some results such as:    generating functions,  Srivastava-Agarwal type generating  functions and   transformational identity  involving the generalized Cauchy polynomials.

The paper is organized as follows:  In Section 2, we state   and prove two theorems on $q$-difference equations.  We give generating functions for   generalized Cauchy polynomials $p_n(x,y,a)$ by using the perspective of $q$-difference equations,  in Section 3.  In Section 4, we derive Srivastava-Agarwal type generating functions involving the generalized Cauchy polynomials. Finally, we obtain a transformational identity involving generating functions for generalized Cauchy polynomials  by the method of homogeneous $q$-difference equations in Section  5.
 
\section{$q$-difference equations}
In this section, we give and prove two theorems  to  be used in the sequel.
\begin{theorem}
\label{prodpos}
Let $f(a,x,y)$ be a three-variable analytic function in a neighborhood of $(a,x,y)=(0,0,0)\in\mathbb{C}^3$. 
 If $f(a,x,y)$  can be expanded in terms of $p_n(x,y,a)$ if and only if 
\be 
\label{aEZA}
x\Big[f(a,x,y)-f(a,x,qy)\Big] 
 =  y \Big[f(a,qx, qy)-f(a,x,qy)\Big]- ay \Big[f(a, qx,q^2y)-f(a,x,q^2y)\Big].
\ee
\end{theorem}
To determine if a given function is an analytic function in several complex variables, we often use the following Hartogs's Theorem. For more information, please refer to
Taylor \cite[p. 28]{Taylor} and Liu \cite[Theorem 1.8]{Liu}.
\begin{lemma} 
  $\big[$Hartogs's Theorem \cite[p.15]{Gunning}$\big]$
\label{lemma1}
If a complex-valued function is holomorphic (analytic) in each variable separately in an open domain $D\subset\mathbb{C}^n$, then it is holomorphic (analytic) in $D.$
\end{lemma}
\begin{lemma}\cite[p. 5 Proposition 1]{Malgrange}
\label{lemma2}
If $f(x_1, x_2, . . . , x_k)$ is analytic at the origin $
(0, 0, . . . , 0)\in\mathbb{C}^k$, then, $f$ can be expanded in an absolutely convergent power series
\be
f(x_1, x_2, . . . , x_k)=\sum_{n_1,n_2,\cdots,n_k=0}^\infty \alpha_{n_1,n_2,\cdots,n_k} x_1^{n_1} x_2^{n_2}\cdots  x_k^{n_k}.
\ee

\end{lemma}
\begin{proof}[Proof of Theorem \ref{prodpos}] From the Hartogs's Theorem and the theory of several complex variables (see Lemmas \ref{lemma1} and \ref{lemma2}), we assume that
\be
\label{a120}
f(a,x,y)=\sum_{k=0}^\infty A_k(a,x)y^k.
\ee
Substituting   (\ref{a120}) into (\ref{aEZA}) yields
\be 
\label{aEA}
 x \sum_{k=0}^\infty (1-q^k)A_k(a,x)y^k  
 =-\sum_{k=0}^\infty (1-aq^k) q^k\Big[A_k(a,x)-A_k(a,qx)\Big] y^{k+1}.
\ee
Comparing coefficients of $y^k,\,k\geq 1$, we readily find that
\be 
 x  (1-q^k)A_k(a,x)  
 =-  (1-aq^{k-1})q^{k-1} \Big[A_{k-1}(a,x)-A_{k-1}(a,qx)\Big]
\ee
which equals to
\bea
\label{aZA}
   A_k(a,x)=   -q^{k-1}\frac{1-aq^{k-1} }{1-q^k} {D}_{q}\Big\{ A_{k-1}(a,x)\Big\}.
\eea
By iteration, we gain
\bea
\label{aZddA}
   A_k(a,x)=  (-1)^kq^{({}^k_2)} \frac{(a;q)_k }{(q;q)_k}{D}_{q}^k\Big\{ A_0(a,x)\Big\}.
\eea
Letting $\displaystyle f(a,x,0)=A_0(a,x)=\sum_{n=0}^\infty \mu_nx^n,$ we have
\bea
\label{dA}
   A_k(a,x)=    (-1)^kq^{({}^k_2)}\frac{(a;q)_k }{(q;q)_k}\sum_{n=0}^\infty \mu_n\frac{(q;q)_n }{(q;q)_{n-k}}x^{n-k}.
\eea
Replacing  (\ref{dA}) in (\ref{a120}), we have:
\bea
f(a,x,y)&=&\sum_{k=0}^\infty  (-1)^kq^{({}^k_2)}\frac{(a;q)_k }{(q;q)_k}\sum_{n=0}^\infty \mu_n\frac{(q;q)_n }{(q;q)_{n-k}}x^{n-k}y^k\cr
&=&\sum_{n=0}^\infty \mu_n\sum_{k=0}^n \qbinomial{n}{k}{q} (-1)^kq^{({}^k_2)}(a;q)_k x^{n-k}y^k. 
\eea
On the other hand, if  $f(a,x,y)$ can be expanded  in term of $ p_n(x,y,a)$, we can verify that  $f(a,x,y)$ satisfies (\ref{aEZA}). The proof of the assertion (\ref{aEZA}) of Theorem \ref{prodpos} is now completed. 
\end{proof}

\begin{theorem}
\label{aprodpos}
Let $f(a,x,y,z)$ be a four-variable analytic function in a neighborhood of $(a,x,y,z)=(0,0,0,0)\in\mathbb{C}^4$. 
\begin{enumerate}
\item   If $f(a,x,y)$  satisfies the $q$-difference equation
\begin{multline}
\label{aaEZA}
 x\Big[f(a,x,y)-f(a,x,qy)\Big] 
 =  y \Big[f(a,qx, qy)-f(a,x,qy)\Big]- ay \Big[f(a, qx,q^2y)-f(a,x,q^2y)\Big]
\end{multline}
then we have:
\bea
\label{aahha}
f(a,x,y)=\widetilde{E}(a,y; {D}_q)\Big\{f(a,x,0) \Big\}.
\eea
\item 
    If $f(a,x,y,z)$  satisfies the $q$-difference equation
\begin{multline}
\label{sEZA}
 (q^{-1}x- y)\Big[f(a, x,y, z)-f(a, x,y,qz)\Big]\\
 =  z \Big[f(a, q^{-1}x,y, qz)-f(a, x, qy, qz)\Big]+ az \Big[f(a, x, qy,q^2z)-f(a,  q^{-1}x,y,q^2z)\Big]
\end{multline}
then we have:
\bea
\label{shha}
f(a, x,y, z)=\widetilde{L}(a,z;\theta_{xy})\Big\{f(a ,x,y, 0) \Big\}.
\eea
  \end{enumerate}
 
\end{theorem}

\begin{corollary}
Let $f(a,b)$ be a two-variable analytic function in a neighborhood of $(a,b)=(0,0)\in\mathbb{C}^2$. 
\\
  If $f(a,b)$  satisfies the $q$-difference equation
\be 
\label{ll}
af(a,b)- b f(qa, qb)
 =(a -b)f(a,qb)
\ee
then we have:
\bea
f(a,b)=R(b{D}_q)\Big\{f(a,0) \Big\}.
\eea

\end{corollary}
 \begin{remark}
For $x=a,\,y=b$ and $z=0$, the relation (\ref{aaEZA}) reduces to (\ref{ll}).\\ 
  For $a=0, x=a, y=0$ and $z=b$, the $q$-difference equation (\ref{sEZA}) reduces to   (\ref{jjZA}).
  \end{remark}
\begin{proof}[Proof of Theorem \ref{aprodpos}]  From the theory of several complex variables \cite{Range}, we begin to solve the $q$-difference equation (\ref{aaEZA}). First we may assume that
\be
\label{aa120}
f(a,x,y)=\sum_{k=0}^\infty A_k(a,x)y^k, 
\ee
Substituting this equation into (\ref{aaEZA}) and comparing coefficients of $y^k,\,k\geq 1$, we readily find that
\be 
 x  (1-q^k)A_k(a,x)  
 =-  (1-aq^{k-1})q^{k-1} \Big[A_{k-1}(a,x)-A_{k-1}(a,qx)\Big]
\ee
which equals to
\bea
\label{aaZA}
   A_k(a,x)=   -q^{k-1}\frac{1-aq^{k-1} }{1-q^k} {D}_{q}\Big\{ A_{k-1}(a,x)\Big\}.
\eea
By iteration, we gain
\bea
\label{aaZddA}
   A_k(a,x)=  (-1)^kq^{({}^k_2)} \frac{(a;q)_k }{(q;q)_k}{D}_{q}^k\Big\{ A_0(a,x)\Big\}.
\eea
Now we return to calculate $A_0(a,x)$. Just taking $y=0$ in (\ref{aa120}), we immediately obtain
$A_0(a,x)$ $=f(a,x,0)$.  The proof of the assertion (\ref{aahha}) of Theorem \ref{aprodpos} is now completed by  substituting (\ref{aaZddA}) back into (\ref{aa120}).\\
Similarly,   we begin to solve the $q$-difference equation (\ref{sEZA}). First we may assume that
\be
\label{120}
f(a,x,y, z)=\sum_{n=0}^\infty B_n(a,x,y)z^n.
\ee
Then  substituting the above equation into (\ref{sEZA}), we have: 
\begin{multline}
(q^{-1}x- y)\sum_{n=0}^\infty (1-q^n)B_n(a, x,y)z^n 
=  \sum_{n=0}^\infty q^n(1-aq^n)[B_n(a ,q^{-1}x,y)-B_n(a, x,qy) 
 ]z^{n+1}
\end{multline}
Comparing coefficients of $z^n,\,n\geq 1$, we readily find that
\be 
(q^{-1}x- y)  (1-q^n)B_n(a,x,y) 
=   q^{n-1}(1-aq^{n-1})[B_{n-1}(a, q^{-1}x,y)-B_{n-1}(a, x,qy) 
].
\ee
After simplification, we get 
\bea
\label{sZA}
  B_n(a,x,y)=    q^{n-1}\frac{1 -aq^{n-1} }{1-q^n} \theta_{xy}\Big\{ B_{n-1}(a,x,y)\Big\}.
\eea
By iteration, we gain
\bea
\label{ZddA}
  B_n(a,x,y)=   \frac{q^{({}^n_2)}(a;q)_n}{(q;q)_n}\theta_{xy}^n\Big\{ B_0(a,x,y)\Big\}.
\eea
Now we return to calculate $A_0(a,x,y)$. Just taking $z=0$ in (\ref{120}), we immediately obtain
$A_0(a,x,y)$ $=f(a,x,y,0)$. The proof of the assertion (\ref{shha}) of Theorem \ref{aprodpos} is now completed by  substituting (\ref{ZddA}) back into (\ref{120}).
\end{proof}

 \section{Generating functions for generalized Cauchy polynomials  }
{\rm The generalized    Cauchy polynomials $p_n(x,y,a)$ \cite{ArjikaSri} are defined as
\be
\label{defR}
p_n(x,y,a)=\sum_{k=0}^n {n \mychoose k}_q\; (-1)^k\;
q^{\binom{k}{2}}\; (a;q)_{k}\; x^{n-k}\; y^k
\ee
}
and their generating function 
\begin{lemma}\cite[Eq. (2.21)]{ArjikaSri}
\label{LM1}
Suppose that $ |xt|<1$, we have:
\be 
\label{dene}
 \sum_{n=0}^\infty p_{n }(x,y,a) \frac{t^n}{(q;q)_n}
=\frac{1}{(xt;q)_\infty}  {}_{1}\Phi_1\left[\begin{array}{rr}a;\\
 \\
0;
 \end{array} 
q; yt  \right].
\ee
\end{lemma}
For $a=0$, in Lemma \ref{LM1}, we get the following
\begin{lemma} \cite{Chen2003}
Suppose that $ |xt|<1$, we have:
\be 
\label{edene}
 \sum_{n=0}^\infty p_{n }(x,y) \frac{t^n}{(q;q)_n}
=\frac{(yt;q)_\infty}{(xt;q)_\infty}.
\ee
\end{lemma}

In this section,   we use the representation    (\ref{defR}) to derive  another generating function  for  generalized Cauchy polynomials by the method of homogeneous $q$-difference equations. 
\begin{theorem}
\label{TA}
Suppose that $ |rx|<1$, we have:
\bea
\label{sums}
\sum_{n=0}^\infty p_n(x,y,a)\frac{(s/r;q)_n\,r^n}{(q;q)_n}=\frac{(sx;q)_\infty}{(rx;q)_\infty} {}_2\Phi_2\left[
\begin{array}{rr}a,s/r;\\\\
sx,0;\end{array} q;ry
\right].
\eea
 \end{theorem}
 \begin{corollary}
 \bea
\label{esums}
\sum_{n=0}^\infty p_n(x,y,a)(-1)^nq^{({}^n_2)}\frac{ \,s^n}{(q;q)_n}= (sx;q)_\infty  {}_1\Phi_2\left[
\begin{array}{rr}a; \\
sx,0;\end{array} q;sy
\right].
\eea
 \end{corollary}
 \begin{remark} 
 For $s=0$ and $r=t$ in Theorem  \ref{TA}, (\ref{sums}) reduces to  (\ref{dene}). For $s=0,\,r=t$  and $a=0$ in Theorem  \ref{TA}, (\ref{sums}) reduces to  (\ref{edene}).
For $r=0$ in Theorem  \ref{TA}, (\ref{sums}) reduces to  (\ref{esums}).
 \end{remark}
 \begin{proof}[Proof of Theorem \ref{TA}]
By denoting the right-hand side of   (\ref{sums}) by $f(a, x, y)$,  we can verify that
$f(a, x, y)$ satisfies (\ref{aEZA}). So, we have
\be
f(a,x,y)=\sum_{n=0}^\infty \mu_n p_n(x,y,a)
\ee
and 
\be
 f(a,x,0) =\sum_{n=0}^\infty \mu_nx^n=\frac{(sx;q)_\infty}{(rx;q)_\infty}=\sum_{n=0}^\infty \frac{(s/r;q)_n\,(rx)^n}{(q;q)_n}.
\ee
 So, $f(a,x,y)$ is equal to the right-hand side of (\ref{sums}).   
 \end{proof}
\begin{theorem} 
\label{orro}
For $k\in\mathbb{N}$ and  $ |xt|<1$, we have:
\be 
\label{ggdene}
 \sum_{n=0}^\infty p_{n+k}(x,y,a) \frac{t^n}{(q;q)_n}
=\frac{x^k}{(xt;q)_\infty}\sum_{n=0}^\infty \frac{(q^{-k}, xt,a;q)_n (yx^{-1}q^k)^n }{(q;q)_n} {}_{1}\Phi_1\left[\begin{array}{rr}aq^{n};\\
 \\
0;
 \end{array}  
q; ytq^n\right].
\ee
   
\end{theorem}
\begin{remark}
For $k=0,$ in Theorem \ref{orro}, (\ref{ggdene}) reduces to (\ref{dene}).
\end{remark}
 \begin{proof}[Proof of Theorem \ref{orro}]
Denoting the right-hand side of equation (\ref{ggdene}) equivalently by 
\be 
\label{vassr}
f(a,  x, y)
 = x^k \sum_{n=0}^\infty \frac{(q^{-k}, xt,a;q)_n (yx^{-1}q^k)^n }{  (q;q)_n}\frac{1  }{(xt q^n;q)_\infty }  {}_{1}\Phi_1\left[\begin{array}{rr}aq^{n};\\
 \\
0;
 \end{array}  
q; ytq^n\right]
\ee
and it is easy to check that (\ref{vassr}) satisfies  (\ref{aaEZA}), so we have:
\be 
f(a,  x, y)
 =  \sum_{n=0}^\infty\mu_n\, p_n(x,y,a).
\ee
Setting $y=0$ in (\ref{vassr}),  it becomes
\be 
f(a,  x, 0)
 =\sum_{n=0}^\infty\mu_n x^n= \frac{ x^k   }{(xt;q)_\infty }=  \sum_{n=0}^\infty x^{n+k} \frac{t^n }{(q,q)_n}= \sum_{n=k}^\infty x^{n } \frac{t^{n-k} }{(q,q)_{n-k}}.
\ee
Hence
\be 
f(a,  x, y)
 =  \widetilde{E}(a,\mu; {D}_{q})\left\{\sum_{n=k}^\infty x^n \frac{t^{n-k} }{(q,q)_{n-k}}\right\}= \sum_{n=k}^\infty \, p_n(x,y,a) \frac{t^{n-k} }{(q,q)_{n-k}}=\sum_{n=0}^\infty \, p_{n+k}(x,y,a) \frac{t^{n } }{(q,q)_{n }},
\ee
which is the left-hand side of  (\ref{ggdene}). 
\end{proof} 

\section{Srivastava-Agarwal type generating functions involving  generalized Cauchy polynomials }

The Hahn polynomials \cite{Hahn049,Hahn49}   (or Al-Salam and Carlitz polynomials \cite{AlSalam}) are given by
\be 
\phi_n^{(a)}(x|q)=\sum_{k=0}^n \qbinomial{n}{k}{q} (a;q)_k x^k.
\ee 
Srivastava and Agarwal deduced the following generating function  (also called
Srivastava-Agarwal type generating functions).
  \begin{lemma} \cite[eq. (3.20)]{SrivastavaAgarwal}  
\be
\label{c1sums}
\sum_{n=0}^\infty \phi_n^{(\alpha)}(x|q) (\lambda;q)_n\frac{  t^n}{(q;q)_n} =  \frac{(\lambda t; q)_\infty }{(t;q)_\infty}   {}_2\Phi_1\left[
\begin{array}{rr} \lambda, \alpha;\\\\
 \lambda  t; \end{array}q; xt   
\right],\quad max\{|t|, |xt|\}<1.
\ee 
 \end{lemma}
 For $\lambda=0,$  we have:
  \begin{lemma} \cite[eq.(1.14)]{JCao}   
\be 
\label{svs}
\sum_{k=0}^\infty  
 \phi_k^{(\alpha)}(x|q)\frac{  t^k}{(q;q)_k}  = \frac{(\alpha xt; q)_\infty }{(xt, t;q)_\infty},\quad max\{|xt|, |t|\}<1.
\ee 
 \end{lemma}
 For more information about Srivastava-Agarwal type generating functions for
Al-Salam-Carlitz polynomials, please refer to \cite{SrivastavaAgarwal,JCao12}.

In this section, we use the representation    (\ref{defR}) to derive Srivastava-Agarwal type generating function  for generalized Cauchy polynomials  by the method of homogeneous $q$-difference equations.
\begin{theorem}
\label{TA1}
For $ M\in\mathbb{N}$, if $\alpha=q^{-M}$ and  $max\{|\lambda t|, |\lambda xt|\}<1$, we have:  
\begin{multline}
\label{1sums}
\sum_{n=0}^\infty \phi_n^{(\alpha)}(x|q)p_n(\lambda, \mu,a)\frac{t^n}{(q;q)_n}\\
= \frac{(\alpha\lambda xt; q)_\infty }{(\lambda xt, \lambda t;q)_\infty}\sum_{k=0}^\infty     \frac{(-1)^{k}q^{({}^{k}_2)}\,(a,\alpha,\lambda t;q)_{k } (\mu x t)^{k }   }{(\alpha\lambda xt, q;q)_{k}}\,    {}_1\Phi_1\left[
\begin{array}{rr} aq^k;\\\\
 0;\end{array}q;    \mu tq^{k}
\right].
\end{multline}
 \end{theorem}
 \begin{remark}
 Setting $a=0, \,\lambda=1$ and $\mu=0$, formula (\ref{1sums}) reduces to (\ref{svs}). For $a=0, \,\lambda=1$ and $\mu=\lambda$, formula (\ref{1sums}) reduces to (\ref{c1sums}).
 \end{remark}
 \begin{proof}[Proof of Theorem \ref{TA1}]  
Denoting the right-hand side of equation (\ref{1sums}) by $H(a,  \lambda, \mu,\alpha,x)$, then we have:
\begin{multline}
\label{ass}
H(a,  \lambda, \mu,\alpha,x)\\
 =\frac{1 }{(\lambda xt ;q)_\infty}\sum_{n=0}^\infty     \frac{(-1)^{n}q^{({}^{n}_2)}\,(a;q)_{n } (\mu  t)^{n }   }{(q;q)_{n}}\, \sum_{k=0}^\infty     \frac{(-1)^{k}q^{({}^{k}_2)}\,(\alpha,aq^n;q)_{k } (\mu x  tq^n)^{k }   }{(q;q)_{k}}\frac{(\alpha\lambda xt q^k; q)_\infty }{(  \lambda tq^k;q)_\infty}.
\end{multline}
We suppose that the operator  ${D}_{q}$ acts upon the variable $\lambda.$  Because equation (\ref{ass}) satisfies (\ref{aaEZA}), we have:
\bea
\label{ss}
H(a,  \lambda, \mu,\alpha,x)= \widetilde{E}(a,\mu; {D}_{q})   \left\{
H(a, \lambda, 0,\alpha,x)\right\} 
&=& \widetilde{E}(a,\mu; {D}_{q})\Bigg\{\frac{(\alpha\lambda xt; q)_\infty }{(\lambda xt, \lambda t;q)_\infty}\Bigg\}\cr
&=&\widetilde{E}(a,\mu; {D}_{q})\Bigg\{\sum_{k=0}^\infty  
 \Phi_k^{(\alpha)}(x|q)\frac{(\lambda t)^k}{(q;q)_k}\Bigg\} 
 \cr
&=&\sum_{k=0}^\infty  
 \Phi_k^{(\alpha)}(x|q)\frac{  t^k}{(q;q)_k}\widetilde{E}(a,\mu; {D}_{q})\{ \lambda ^k\} \nonumber
\eea
which is the left-hand side of  (\ref{ass}). The proof is complete. 
\end{proof} 
\section{A transformational identity involving generating functions for generalized Cauchy polynomials }
In this section we deduce the following transformational identity involving generating functions for generalized Cauchy polynomials  by the method of homogeneous
$q$-difference equation.
\begin{theorem}
\label{asss}
Let $A(k)$ and $B(k)$  satisfy
\be
\label{sam}
\sum_{k=0}^\infty A(k) x^k= \sum_{k=0}^\infty B(k) \frac{1}{(xtq^k;q)_\infty}
\ee
and we have
\be
\label{samm}
\sum_{k=0}^\infty A(k) p_k(x,y,a)= \sum_{k=0}^\infty B(k) \frac{1}{(xtq^k;q)_\infty}{}_1\Phi_1\left[\begin{array}{rr}a;\\\\
0;\end{array}  q;ytq^k\right]
\ee
supposing that (\ref{sam}) and  (\ref{samm}) are convergent.
\end{theorem}
\begin{proof} 
We denote the right-hand side of (\ref{samm})  by $f(a,x,y)$  and we   can check that  $f(a,x,y)$ satisfies  (\ref{aaEZA}). We then obtain 
\be 
\label{vass}
f(a,  x, y)
 =  \sum_{k=0}^\infty\mu_k p_k(x,y,a)
\ee
and 
\bea 
f(a,  x, 0)
 =\sum_{k=0}^\infty\mu_k x^k&=& \sum_{k=0}^\infty B(k) \frac{1}{(xtq^k;q)_\infty}\, (\mbox{by (\ref{sam})})\cr
 &=&\sum_{k=0}^\infty   A(k) x^k .
\eea
Hence
\be 
f(a,  x, y)
 =  \sum_{k=0}^\infty A(k)\, p_k(x,y,a),
\ee
which is the left-hand side of  (\ref{samm}). The proof of Theorem \ref{asss} is thus completed.
\end{proof}


\begin{thebibliography}{99}
\bibitem{AlSalam}W. A. AL-Salam  and L Carlitz, {\it Some orthogonal $q$-polynomials}, Math. Nachr. {\bf 30 } (1965), pp. 47--61.

\bibitem{Andrews}  G. E. Andrews,     {\it  The theory of partitions}, Cambridge Univ. Press, 1985.

\bibitem{JCao12} J. Cao, {\it Generalizations of certain Carlitz's trilinear and Srivastava-Agarwal type generating functions},  J. Math. Anal. Appl. {\bf 396} (2012), pp.  351--362.

\bibitem{JCao16} J. Cao, {\it Homogeneous $q$-difference equations and generating
functions for $q$-hypergeometric polynomials}, 
Ramanujan J. {\bf 40} (2016), pp. 177--192.

\bibitem{JCao}J. Cao, {\it Alternative proofs of generating functions for Hahn polynomials and some applications}, Infin. Dimens. Anal. Quantum Probab. Relat. Top. {\bf 14 } (2011), pp. 571--590.

\bibitem{Chen2003}   W. Y. C. Chen,  A. M. Fu and    B. Zhang,    {\it The homogeneous $q$-difference operator}, Adv. App. Math. {\bf 31} (2003), pp.   659--668.

\bibitem{Liu97}W. Y. C.  Chen and     Z.-G. Liu,    {\it Parameter augmenting for basic hypergeometric series II} J. Combin. Theory, Ser. A, {\bf 80}, (1997)  pp. 175--195.

\bibitem{Liu98}W. Y. C.   Chen and     Z.-G. Liu,   {\it  Parameter augmenting for basic hypergeometric series  I}, in  Mathematical Essays in Honor of Gian-Carlo Rota, B. E. Sagan and R.P. Stanley, eds., Birk\"auser, Basel, (1998), pp. 111--129.

\bibitem{GasparRahman} G.  Gasper  and  M. Rahman,  {\it Basic Hypergeometric Series}, 2nd edn. Cambridge University Press, Cambridge, 2004.

\bibitem{Goldman} J. Golman and  G.-C. Rota, {\it On the foundations of combinatorial theory, IV: Finite vector spaces and Eulerien generating functions}, Sut. Appl. Math. {\bf 49}  (1970), pp. 239--258. 

\bibitem{Gunning} R. Gunning, {\it Introduction to Holomorphic Functions of Several Variables}. In: Function theory, vol. {\bf 1}, Wadsworth and Brooks/Cole, Belmont, 1990.

\bibitem{Hahn049} W. Hahn, {\it Uber Orthogonalpolynome, die $q$-Differenzengleichungen}, Math. Nuchr. {\bf 2} (1949), 434.

\bibitem{Hahn49} W. Hahn,   {\it Beitrage zur Theorie der Heineschen Reihen; Die 24 Integrale der hypergeometrischen $q$-Differenzengleichung; Das $q$-Analogon der 
Laplace-Transformation}, Math. Nachr. {\bf 2} (1949), pp. 340--379.

\bibitem{Ismail1981}  E. C. Ihrig and      M. E. H. Ismail,    { A $q$-umbral calculus}, J. Math. Anal. Appl. {\bf 84} (1981), pp. 178--207. 

\bibitem{Johnson}   W. P. Johnson,     {\it $q$-Extensions of identities of Abel-Rothe type} Discrete Math. {\bf 159}(1995), pp. 161--177.  

\bibitem{Koekock}R. Koekock  and R. F.  Swarttouw,   {\it The Askey-scheme of hypergeometric orthogonal polynomials and its $q$-analogue report}, Delft University of Technology, 1998.

\bibitem{Liu10} Z.-G. Liu, {\it Two $q$-difference equations and $q$-operator identities}, J. Differ. Equ. Appl. {\bf 16 } (2010), pp. 1293--1307.

\bibitem{Liu11}Z.-G. Liu, {\it An extension of the non-terminating \textcolor{red}{${}_6\phi_5$} summation and the Askey-Wilson polynomials}, J. Differ. Equ. Appl. {\bf 17}  (2011), pp. 1401--1411.

\bibitem{Liu}Z.-G. Liu,  {\it On the $q$-partial Differential Equations and $q$-series}. In: The legacy of Srinivasa Ramanujan, Ramanujan Mathematical Society Lecture Note Series, Mysore   {\bf 20} (2013), pp. 213--250.

\bibitem{Malgrange}B. Malgrange,   {\it Lectures on the Theory of Functions of Several Complex Variables}, Springer, Berlin, 1984.

\bibitem{Range} R. M. Range, {\it Complex analysis: A brief tour into higher dimensions}, Amer. Math. Mon. {\bf 110} (2003), pp. 89--108.

\bibitem{Roman1982}   S. Roman,   {\it The theory of the umbral calculus I},  J. Math. Anal. Appl. {\bf 87} (1982),  pp. 58--115.

\bibitem{Saadsukhi} H. L.    Saad and    A. A. Sukhi,  {\it Another homogeneous $q$-difference operator}, Appl. Math. Comput. {\bf 215} (2010), 4332--4339. 

\bibitem{Saad} H. L.    Saad and    A. A. Sukhi,   {\it  The $q$-Exponential Operator}, Appl. Math. Sci. {\bf 7} (2005), 6369--6380. 

\bibitem{SLATER} L. J. Slatter,  {\it Generalized Hypergeometric Functions}, Cambridge Univ. Press, Cambridge, London and New York,  1966. 
  
\bibitem{SrivastaKarlsson} H. M. Srivastava and P. W. Karlsson,  
{\it Multiple Gaussian Hypergeometric Series}, Halsted Press (Ellis Horwood Limited, Chichester);  John Wiley and Sons, New York, Chichester, Brisbane and Toronto, 1985. 
 
\bibitem{SrivastaAbdlhusein} H. M. Srivastava and   M. A. Abdlhusein,   {\it New forms of the Cauchy operator and some of their  applications},  Russian J. Math. Phys. {\bf 23} (2016), 124--134.

\bibitem{SrivastavaAgarwal} H. M. Srivastava and A.K. Agarwal, {\it Generating functions for a class of $q$-polynomials}, Ann. Mat. Pure Appl. (Ser. 4) {\bf 154} (1989), pp. 99--109.

\bibitem{ArjikaSri}  H. M. Srivastava,  S. Arjika and   A. Sherif Kelil,   {\it Some homogeneous $q$-difference operators and the associated generalized Hahn polynomials}, Appl. Set-Valued Anal. Optim. {\bf 1} (2019), 187--201. 

\bibitem{Taylor}J. Taylor,  {\it  Several complex variables with connections to algebraic geometry and lie groups}, Graduate Studies in Mathematics, American Mathematical Society, Providence, vol. {\bf 46}  2002.  

\end{thebibliography}
\end{document}